\documentclass{amsart}
\usepackage{amssymb}

\usepackage{amsmath}
\usepackage{amscd}
\usepackage{thmdefs}

\input{tcilatex}
\theoremstyle{definition}
\theoremstyle{remark}
\numberwithin{equation}{section}

\input tcilatex

\begin{document}
\title[Compressed Random Variables]{Compressed Random Variables in the Graph $W^{*}$-Probability Spaces}
\author{Ilwoo Cho}
\address{Univ. of Iowa, Dep. of Math, Iowa City, IA, U. S. A.\\
}
\email{ilcho@math.uiowa.edu}
\keywords{Graph $W^{*}$-Probability Spaces over Diagonal Subalgebras,
Diagonal-Compressed Random Variables, Off-Diagonal Compressed Random
Variables. Compressed Random Variables.}
\maketitle

\begin{abstract}
In [16] and [17], we observed the amalgamated free probability theory on the
graph $W^{*}$-probability space $\left( W^{*}(G),E\right) $ over the
diagonal subalgebra $D_{G}.$ In [18], we consider the diagonal compressed
random variables in $\left( W^{*}(G),E\right) $ and observed the amalgamated
freeness in the $W^{*}$-probability space $\left( W^{*}(G),E\right) $ over
the $D_{G}.$ In particular, for $V_{N}=\{v_{0}\}\subset V(G),$ we could
construct the (scalar-valued) tracial $W^{*}$-probability space, so-called
the vertex compressed graph $W^{*}$-probability space. In this paper, we
will consider the off-diagonal compressed random variables in $\left(
W^{*}(G),E\right) .$ After fixing $v_{1}\neq v_{2}$ in $V(G),$ we define the 
$(v_{1},v_{2})$-off-diagonal compressed random variable of $a\in \left(
W^{*}(G),E\right) ,$ $L_{v_{1}}aL_{v_{2}}.$ We will consider the free
probability data on such elements in $\left( W^{*}(G),E\right) .$ Also, we
will consider the compressed random variables $PaP$ for the $D_{G}$-valued
random variable $a\in \left( W^{*}(G),E\right) ,$ compressed by the
projection $P=L_{v_{1}}+L_{v_{2}}+...+L_{v_{N}}\in D_{G}.$ To observe the
free probability data for such compressed random variable, we use the
diagonal-compression and off-diagonal compression. We can figure out that
only the diagonal-compression affects the compressed free probability on $%
W^{*}(G).$ In fact, $D_{G}$-valued moment series and R-transform of the
compressed random variable of $a$ are same as those of diagonal compressed
random of it.
\end{abstract}

\strut

In [16], we constructed the graph $W^{*}$-probability spaces. The graph $%
W^{*}$-probability theory is one of the good example of Speicher's
combinatorial free probability theory with amalgamation. In [16], we
observed how to compute the moment and cumulant of an arbitrary random
variables in the graph $W^{*}$-probability space and the freeness on it with
respect to the given conditional expectation. Also, in [17], we consider
certain special random variables of the graph $W^{*}$-probability space, for
example, semicircular elements, even elements and R-diagonal elements. This
shows that the graph $W^{*}$-probability spaces contain the rich free
probabilistic objects. \bigskip Roughly speaking, graph $W^{*}$-algebras are 
$W^{*}$-topology closed version of free semigroupoid algebras defined and
observed by Kribs and Power in [10].

\strut

Throughout this paper, let $G$ be a countable directed graph and let $%
\mathbb{F}^{+}(G)$ be the free semigroupoid of $G,$ in the sense of Kribs
and Power. i.e., it is a collection of all vertices of the graph $G$ as
units and all admissible finite paths, under the admissibility. As a set,
the free semigroupoid $\mathbb{F}^{+}(G)$ can be decomposed by

\strut

\begin{center}
$\mathbb{F}^{+}(G)=V(G)\cup FP(G),$
\end{center}

\strut

where $V(G)$ is the vertex set of the graph $G$ and $FP(G)$ is the set of
all admissible finite paths. Trivially the edge set $E(G)$ of the graph $G$
is properly contained in $FP(G),$ since all edges of the graph can be
regarded as finite paths with their length $1.$ We define a graph $W^{*}$%
-algebra of $G$ by

\strut

\begin{center}
$W^{*}(G)\overset{def}{=}\overline{%
\mathbb{C}[\{L_{w},L_{w}^{*}:w\in
\mathbb{F}^{+}(G)\}]}^{w},$
\end{center}

\strut

where $L_{w}$ and $L_{w}^{\ast }$ are creation operators and annihilation
operators on the generalized Fock space $H_{G}=l^{2}\left( \mathbb{F}%
^{+}(G)\right) $ induced by the given graph $G,$ respectively. Notice that
the creation operators induced by vertices are projections and the creation
operators induced by finite paths are partial isometries. We can define the $%
W^{\ast }$-subalgebra $D_{G}$ of $W^{\ast }(G),$ which is called the
diagonal subalgebra by

\strut

\begin{center}
$D_{G}\overset{def}{=}\overline{\mathbb{C}[\{L_{v}:v\in V(G)\}]}^{w}.$
\end{center}

\strut

Then each element $a$ in the graph $W^{*}$-algebra $W^{*}(G)$ is expressed by

\strut

\begin{center}
$a=\underset{w\in \mathbb{F}^{+}(G:a),\,u_{w}\in \{1,*\}}{\sum }%
p_{w}L_{w}^{u_{w}},$ \ for $p_{w}\in \mathbb{C},$
\end{center}

\strut

where $\mathbb{F}^{+}(G:a)$ is a support of the element $a$, as a subset of
the free semigroupoid $\mathbb{F}^{+}(G).$ The above expression of the
random variable $a$ is said to be the Fourier expansion of $a.$ Since $%
\mathbb{F}^{+}(G)$ is decomposed by the disjoint subsets $V(G)$ and $FP(G),$
the support $\mathbb{F}^{+}(G:a)$ of $a$ is also decomposed by the following
disjoint subsets,

\strut

\begin{center}
$V(G:a)=\mathbb{F}^{+}(G:a)\cap V(G)$
\end{center}

and

\begin{center}
$FP(G:a)=\mathbb{F}^{+}(G:a)\cap FP(G).$
\end{center}

\strut

Thus the operator $a$ can be re-expressed by

\strut

\begin{center}
$a=\underset{v\in V(G:a)}{\sum }p_{v}L_{v}+\underset{w\in FP(G:a),\,u_{w}\in
\{1,*\}}{\sum }p_{w}L_{w}^{u_{w}}.$
\end{center}

\strut

Notice that if $V(G:a)\neq \emptyset ,$ then $\underset{v\in V(G:a)}{\sum }%
p_{v}L_{v}$ is contained in the diagonal subalgebra $D_{G}.$ Thus we have
the canonical conditional expectation $E:W^{*}(G)\rightarrow D_{G},$ defined
by

\strut

\begin{center}
$E\left( a\right) =\underset{v\in V(G:a)}{\sum }p_{v}L_{v},$
\end{center}

\strut

for all $a=\underset{w\in \mathbb{F}^{+}(G:a),\,u_{w}\in \{1,*\}}{\sum }%
p_{w}L_{w}^{u_{w}}$ \ in $W^{*}(G).$ Then the algebraic pair $\left(
W^{*}(G),E\right) $ is a $W^{*}$-probability space with amalgamation over $%
D_{G}$ (See [16]). It is easy to check that the conditional expectation $E$
is faithful in the sense that if $E(a^{*}a)=0_{D_{G}},$ for $a\in W^{*}(G),$
then $a=0_{D_{G}}.$

\strut

For the fixed operator $a\in W^{*}(G),$ the support $\mathbb{F}^{+}(G:a)$ of
the operator $a$ is again decomposed by

\strut

\begin{center}
$\mathbb{F}^{+}(G:a)=V(G:a)\cup FP_{*}(G:a)\cup FP_{*}^{c}(G:a),$
\end{center}

\strut

with the decomposition of $FP(G:a),$

$\strut $

\begin{center}
$FP(G:a)=FP_{*}(G:a)\cup FP_{*}^{c}(G:a),$
\end{center}

where

\strut

\begin{center}
$FP_{*}(G:a)=\{w\in FP(G:a):$both $L_{w}$ and $L_{w}^{*}$ are summands of $%
a\}$
\end{center}

and

\begin{center}
$FP_{*}(G:a)=FP(G:a)\,\,\setminus \,\,FP_{*}(G:a).$
\end{center}

\strut

The above new expression plays a key role to find the $D_{G}$-valued moments
of the random variable $a.$ In fact, the summands $p_{v}L_{v}$'s and $%
p_{w}L_{w}+p_{w^{t}}L_{w}^{*},$ for $v\in V(G:a)$ and $w\in FP_{*}(G:a)$ act
for the computation of $D_{G}$-valued moments of $a.$ By using the above
partition of the support of a random variable, we can compute the $D_{G}$%
-valued moments and $D_{G}$-valued cumulants of it via the lattice path
model $LP_{n}$ and the lattice path model $LP_{n}^{*}$ satisfying the $*$%
-axis-property. At a first glance, the computations of $D_{G}$-valued
moments and cumulants look so abstract and hence it looks useless. However,
these computations, in particular the computation of $D_{G}$-valued
cumulants, provides us how to figure out the $D_{G}$-freeness of random
variables by making us compute the mixed cumulants. As applications, in the
final chapter, we can compute the moment and cumulant of the operator that
is the sum of $N$-free semicircular elements with their covariance $2.$

\strut \strut \strut

Based on the $D_{G}$-cumulant computation, we can characterize the $D_{G}$%
-freeness of generators of $W^{*}(G),$ by the so-called diagram-distinctness
on the graph $G.$ i.e., the random variables $L_{w_{1}}$ and $L_{w_{2}}$ are
free over $D_{G}$ if and only if $w_{1}$ and $w_{2}$ are diagram-distinct
the sense that $w_{1}$ and $w_{2}$ have different diagrams on the graph $G.$
Also, we could find the necessary condition for the $D_{G}$-freeness of two
arbitrary random variables $a$ and $b.$ i.e., if the supports $\mathbb{F}%
^{+}(G:a)$ and $\mathbb{F}^{+}(G:b)$ are diagram-distinct, in the sense that 
$w_{1}$ and $w_{2}$ are diagram distinct for all pairs $(w_{1},w_{2})$ $\in $
$\mathbb{F}^{+}(G:a)$ $\times $ $\mathbb{F}^{+}(G:b),$ then the random
variables $a$ and $b$ are free over $D_{G}.$

\strut \strut

In [17], we considered some special $D_{G}$-valued random variables in a
graph $W^{*}$-probability space $\left( W^{*}(G),E\right) .$ The those
random variables are the basic objects to study Free Probability Theory. We
can conclude that

\strut

(i) \ \ if $l$ is a loop, then $L_{l}+L_{l}^{*}$ is $D_{G}$-semicircular.

\strut

(ii) \ if $w$ is a finite path, then $L_{w}+L_{w}^{*}$ is $D_{G}$-even.

\strut

(iii) if $w$ is a finite path, then $L_{w}$ and $L_{w}^{*}$ are $D_{G}$%
-valued R-diagonal.

\strut

In [18], we observed the diagonal compressed random variables in the graph $%
W^{*}$-probability space $\left( W^{*}(G),E\right) .$ Let $%
v_{1},...,v_{N}\in V(G)$ and let $a$ be a $D_{G}$-valued random variable in $%
\left( W^{*}(G),E\right) .$ Define the diagonal compressed random variable
of $a$ by $V=\{v_{1},...,v_{N}\}$ by

\strut

\begin{center}
$C_{V}(a)=L_{v_{1}}aL_{v_{1}}+...+L_{v_{N}}aL_{v_{N}}.$
\end{center}

\strut

Notice that if $v\in V(G),$ then $L_{v}aL_{v}$ is the compressed random
variable by $L_{v}$ and the compressed random variable has its support
contained in $\{v\}$ $\cup $ $loop_{v}(G),$ where $loop_{v}(G)$ $=$ $\{l\in
loop(G):$ $l=vlv\}.$

\strut

The main purpose of this paper is to show that the compressed random
variable $P_{V}aP_{V}$ of a projection $P_{V}=\sum_{j=1}^{N}L_{v_{j}},$
where $V=\{v_{1},...,v_{N}\},$ has the same free probabilistic information
with the diagonal compressed random variable $C_{V}(a).$ i.e., the
compressed random variables $P_{V}aP_{V}$ and $C_{V}(a)$ have the same $%
D_{G} $-valued moments and cumulants.

\strut \strut

\strut \strut

\strut \strut

\section{Graph $W^{*}$-Probability Theory}

\strut

\strut

Let $G$ be a countable directed graph and let $\Bbb{F}^{+}(G)$ be the free
semigroupoid of $G.$ i.e., the set $\mathbb{F}^{+}(G)$ is the collection of
all vertices as units and all admissible finite paths of $G.$ Let $w$ be a
finite path with its source $s(w)=x$ and its range $r(w)=y,$ where $x,y\in
V(G).$ Then sometimes we will denote $w$ by $w=xwy$ to express the source
and the range of $w.$ We can define the graph Hilbert space $H_{G}$ by the
Hilbert space $l^{2}\left( \mathbb{F}^{+}(G)\right) $ generated by the
elements in the free semigroupoid $\mathbb{F}^{+}(G).$ i.e., this Hilbert
space has its Hilbert basis $\mathcal{B}=\{\xi _{w}:w\in \mathbb{F}%
^{+}(G)\}. $ Suppose that $w=e_{1}...e_{k}\in FP(G)$ is a finite path with $%
e_{1},...,e_{k}\in E(G).$ Then we can regard $\xi _{w}$ as $\xi
_{e_{1}}\otimes ...\otimes \xi _{e_{k}}.$ So, in [10], Kribs and Power
called this graph Hilbert space the generalized Fock space. Throughout this
paper, we will call $H_{G}$ the graph Hilbert space to emphasize that this
Hilbert space is induced by the graph.

\strut

Define the creation operator $L_{w},$ for $w\in \mathbb{F}^{+}(G),$ by the
multiplication operator by $\xi _{w}$ on $H_{G}.$ Then the creation operator 
$L$ on $H_{G}$ satisfies that

\strut

(i) \ $L_{w}=L_{xwy}=L_{x}L_{w}L_{y},$ for $w=xwy$ with $x,y\in V(G).$

\strut

(ii) $L_{w_{1}}L_{w_{2}}=\left\{ 
\begin{array}{lll}
L_{w_{1}w_{2}} &  & \text{if }w_{1}w_{2}\in \mathbb{F}^{+}(G) \\ 
&  &  \\ 
0 &  & \text{if }w_{1}w_{2}\notin \mathbb{F}^{+}(G),
\end{array}
\right. $

\strut

\ \ \ \ for all $w_{1},w_{2}\in \mathbb{F}^{+}(G).$

\strut

Now, define the annihilation operator $L_{w}^{*},$ for $w\in \mathbb{F}%
^{+}(G)$ by

\strut

\begin{center}
$L_{w}^{\ast }\xi _{w^{\prime }}\overset{def}{=}\left\{ 
\begin{array}{lll}
\xi _{h} &  & \text{if }w^{\prime }=wh\in \mathbb{F}^{+}(G)\xi \\ 
&  &  \\ 
0 &  & \text{otherwise.}
\end{array}
\right. $
\end{center}

\strut

The above definition is gotten by the following observation ;

\strut

\begin{center}
$
\begin{array}{ll}
<L_{w}\xi _{h},\xi _{wh}>\, & =\,<\xi _{wh},\xi _{wh}>\, \\ 
& =\,1=\,<\xi _{h},\xi _{h}> \\ 
& =\,<\xi _{h},L_{w}^{*}\xi _{wh}>,
\end{array}
\,$
\end{center}

\strut

where $<,>$ is the inner product on the graph Hilbert space $H_{G}.$ Of
course, in the above formula we need the admissibility of $w$ and $h$ in $%
\mathbb{F}^{+}(G).$ However, even though $w$ and $h$ are not admissible
(i.e., $wh\notin \mathbb{F}^{+}(G)$), by the definition of $L_{w}^{\ast },$
we have that

\strut

\begin{center}
$
\begin{array}{ll}
<L_{w}\xi _{h},\xi _{h}> & =\,<0,\xi _{h}> \\ 
& =0=\,<\xi _{h},0> \\ 
& =\,<\xi _{h},L_{w}^{*}\xi _{h}>.
\end{array}
\,\,$
\end{center}

\strut

Notice that the creation operator $L$ and the annihilation operator $L^{*}$
satisfy that

\strut

(1.1) \ \ \ $L_{w}^{*}L_{w}=L_{y}$ \ \ and \ \ $L_{w}L_{w}^{*}=L_{x},$ \ for
all \ $w=xwy\in \mathbb{F}^{+}(G),$

\strut

\textbf{under the weak topology}, where $x,y\in V(G).$ Remark that if we
consider the von Neumann algebra $W^{*}(\{L_{w}\})$ generated by $L_{w}$ and 
$L_{w}^{*}$ in $B(H_{G}),$ then the projections $L_{y}$ and $L_{x}$ are
Murray-von Neumann equivalent, because there exists a partial isometry $%
L_{w} $ satisfying the relation (1.1). Indeed, if $w=xwy$ in $\mathbb{F}%
^{+}(G), $ with $x,y\in V(G),$ then under the weak topology we have that

\strut

(1,2) \ \ \ $L_{w}L_{w}^{*}L_{w}=L_{w}$ \ \ and \ \ $%
L_{w}^{*}L_{w}L_{w}^{*}=L_{w}^{*}.$

\strut

So, the creation operator $L_{w}$ is a partial isometry in $W^{*}(\{L_{w}\})$
in $B(H_{G}).$ Assume now that $v\in V(G).$ Then we can regard $v$ as $%
v=vvv. $ So,

\strut

(1.3) $\ \ \ \ \ \ \ \ \ L_{v}^{*}L_{v}=L_{v}=L_{v}L_{v}^{*}=L_{v}^{*}.$

\strut

This relation shows that $L_{v}$ is a projection in $B(H_{G})$ for all $v\in
V(G).$

\strut

Define the \textbf{graph }$W^{*}$\textbf{-algebra} $W^{*}(G)$ by

\strut

\begin{center}
$W^{*}(G)\overset{def}{=}\overline{%
\mathbb{C}[\{L_{w},L_{w}^{*}:w\in
\mathbb{F}^{+}(G)\}]}^{w}.$
\end{center}

\strut

Then all generators are either partial isometries or projections, by (1.2)
and (1.3). So, this graph $W^{\ast }$-algebra contains a rich structure, as
a von Neumann algebra. (This construction can be the generalization of that
of group von Neumann algebra.) Naturally, we can define a von Neumann
subalgebra $D_{G}\subset W^{\ast }(G)$ generated by all projections $L_{v},$ 
$v\in V(G).$ i.e.

\strut

\begin{center}
$D_{G}\overset{def}{=}W^{*}\left( \{L_{v}:v\in V(G)\}\right) .$
\end{center}

\strut

We call this subalgebra the \textbf{diagonal subalgebra} of $W^{*}(G).$
Notice that $D_{G}=\Delta _{\left| G\right| }\subset M_{\left| G\right| }(%
\mathbb{C}),$ where $\Delta _{\left| G\right| }$ is the subalgebra of $%
M_{\left| G\right| }(\mathbb{C})$ generated by all diagonal matrices. Also,
notice that $1_{D_{G}}=\underset{v\in V(G)}{\sum }L_{v}=1_{W^{*}(G)}.$

\strut

If $a\in W^{*}(G)$ is an operator, then it has the following decomposition
which is called the Fourier expansion of $a$ ;

\strut

(1.4) $\ \ \ \ \ \ \ \ \ \ \ a=\underset{w\in \mathbb{F}^{+}(G:a),\,u_{w}\in
\{1,*\}}{\sum }p_{w}L_{w}^{u_{w}},$

\strut

where $p_{w}\in C$ and $\mathbb{F}^{+}(G:a)$ is the support of $a$ defined by

\strut

\begin{center}
$\mathbb{F}^{+}(G:a)=\{w\in \mathbb{F}^{+}(G):p_{w}\neq 0\}.$
\end{center}

\strut

Remark that the free semigroupoid $\mathbb{F}^{+}(G)$ has its partition $%
\{V(G),FP(G)\},$ as a set. i.e.,

\strut

\begin{center}
$\mathbb{F}^{+}(G)=V(G)\cup FP(G)$ \ \ and \ \ $V(G)\cap FP(G)=\emptyset .$
\end{center}

\strut

So, the support of $a$ is also partitioned by

\strut

\begin{center}
$\mathbb{F}^{+}(G:a)=V(G:a)\cup FP(G:a),$
\end{center}

\strut where

\begin{center}
$V(G:a)\overset{def}{=}V(G)\cap \mathbb{F}^{+}(G:a)$
\end{center}

and

\begin{center}
$FP(G:a)\overset{def}{=}FP(G)\cap \mathbb{F}^{+}(G:a).$
\end{center}

\strut

So, the above Fourier expansion (1.4) of the random variable $a$ can be
re-expressed by

\strut

(1.5) $\ \ \ \ \ \ a=\underset{v\in V(G:a)}{\sum }p_{v}L_{v}+\underset{w\in
FP(G:a),\,u_{w}\in \{1,*\}}{\sum }p_{w}L_{w}^{u_{w}}.$

\strut

We can easily see that if $V(G:a)\neq \emptyset ,$ then $\underset{v\in
V(G:a)}{\sum }p_{v}L_{v}$ is contained in the diagonal subalgebra $D_{G}.$
Also, if $V(G:a)=\emptyset ,$ then $\underset{v\in V(G:a)}{\sum }%
p_{v}L_{v}=0_{D_{G}}.$ So, we can define the following canonical conditional
expectation $E:W^{*}(G)\rightarrow D_{G}$ by

\strut

(1.6) \ \ \ $E(a)=E\left( \underset{w\in \mathbb{F}^{+}(G:a),\,u_{w}\in
\{1,*\}}{\sum }p_{w}L_{w}^{u_{w}}\right) \overset{def}{=}\underset{v\in
V(G:a)}{\sum }p_{v}L_{v},$

\strut

for all $a\in W^{*}(G).$ Indeed, $E$ is a well-determined conditional
expectation.

\strut \strut \strut \strut

\begin{definition}
Let $G$ be a countable directed graph and let $W^{*}(G)$ be the graph $W^{*}$%
-algebra induced by $G.$ Let $E:W^{*}(G)\rightarrow D_{G}$ be the
conditional expectation defined above. Then we say that the algebraic pair $%
\left( W^{*}(G),E\right) $ is the graph $W^{*}$-probability space over the
diagonal subalgebra $D_{G}$. By the very definition, it is one of the $W^{*}$%
-probability space with amalgamation over $D_{G}.$ All elements in $\left(
W^{*}(G),E\right) $ are called $D_{G}$-valued random variables.
\end{definition}

\strut

We have a graph $W^{*}$-probability space $\left( W^{*}(G),E\right) $ over
its diagonal subalgebra $D_{G}.$ We will define the following free
probability data of $D_{G}$-valued random variables.

\strut

\begin{definition}
Let $W^{*}(G)$ be the graph $W^{*}$-algebra induced by $G$ and let $a\in
W^{*}(G).$ Define the $n$-th ($D_{G}$-valued) moment of $a$ by

\strut 

$\ \ \ \ \ E\left( d_{1}ad_{2}a...d_{n}a\right) ,$ for all $n\in \mathbb{N}$,

\strut 

where $d_{1},...,d_{n}\in D_{G}$. Also, define the $n$-th ($D_{G}$-valued)
cumulant of $a$ by

\strut 

$\ \ \ \ \ k_{n}(d_{1}a,d_{2}a,...,d_{n}a)=C^{(n)}\left( d_{1}a\otimes
d_{2}a\otimes ...\otimes d_{n}a\right) ,$

\strut 

for all $n\in \mathbb{N},$ and for $d_{1},...,d_{n}\in D_{G},$ where $%
\widehat{C}=(C^{(n)})_{n=1}^{\infty }\in I^{c}\left( W^{*}(G),D_{G}\right) $
is the cumulant multiplicative bimodule map induced by the conditional
expectation $E,$ in the sense of Speicher. We define the $n$-th trivial
moment of $a$ and the $n$-th trivial cumulant of $a$ by

\strut 

$\ \ \ \ \ E(a^{n})$ $\ \ $and $\ \ k_{n}\left( \underset{n-times}{%
\underbrace{a,a,...,a}}\right) =C^{(n)}\left( a\otimes a\otimes ...\otimes
a\right) ,$

\strut 

respectively, for all $n\in \mathbb{N}.$
\end{definition}

\strut

To compute the $D_{G}$-valued moments and cumulants of the $D_{G}$-valued
random variable $a,$ we need to introduce the following new definition ;

\strut

\begin{definition}
Let $\left( W^{*}(G),E\right) $ be a graph $W^{*}$-probability space over $%
D_{G}$ and let $a\in \left( W^{*}(G),E\right) $ be a random variable. Define
the subset $FP_{*}(G:a)$ in $FP(G:a)$ \ by

\strut 

$\ \ \ FP_{*}\left( G:a\right) \overset{def}{=}\{w\in \mathbb{F}^{+}(G:a):$%
both $L_{w}$ and $L_{w}^{*}$ are summands of $a\}.$

\strut 

And let $FP_{*}^{c}(G:a)\overset{def}{=}FP(G:a)\,\setminus \,FP_{*}(G:a).$
\end{definition}

\strut \strut \strut

We already observed that if $a\in \left( W^{*}(G),E\right) $ is a $D_{G}$%
-valued random variable, then $a$ has its Fourier expansion $a_{d}+a_{0},$
where

\strut

\begin{center}
$a_{d}=\underset{v\in V(G:a)}{\sum }p_{v}L_{v}$
\end{center}

and

\begin{center}
$a_{0}=\underset{w\in FP(G:a),\,u_{w}\in \{1,*\}}{\sum }p_{w}L_{w}^{u_{w}}.$
\end{center}

\strut

By the previous definition, the set $FP(G:a)$ is partitioned by

\strut

\begin{center}
$FP(G:a)=FP_{*}(G:a)\cup FP_{*}^{c}(G:a),$
\end{center}

\strut

for the fixed random variable $a$ in $\left( W^{*}(G),E\right) .$ So, the
summand $a_{0},$ in the Fourier expansion of $a=a_{d}+a_{0},$ has the
following decomposition ;

\strut

\begin{center}
$a_{0}=a_{(*)}+a_{(non-*)},$
\end{center}

\strut where\strut

\begin{center}
$a_{(*)}=\underset{l\in FP_{*}(G:a)}{\sum }\left(
p_{l}L_{l}+p_{l^{t}}L_{l}^{*}\right) $
\end{center}

and

\begin{center}
$a_{(non-*)}=\underset{w\in FP_{*}^{c}(G:a),\,u_{w}\in \{1,*\}}{\sum }%
p_{w}L_{w}^{u_{w}},$
\end{center}

\strut

where $p_{l^{t}}$ is the coefficient of $L_{l}^{*}$ depending on $l\in
FP_{*}(G:a).$

\strut \strut \strut

\strut \strut

\strut \strut

\subsection{$D_{G}$-Moments and $D_{G}$-Cumulants of Random Variables}

\strut

\strut

\strut

Throughout this chapter, let $G$ be a countable directed graph and let $%
\left( W^{*}(G),E\right) $ be the graph $W^{*}$-probability space over its
diagonal subalgebra $D_{G}.$ In this chapter, we will compute the $D_{G}$%
-valued moments and the $D_{G}$-valued cumulants of arbitrary random variable

$\strut $

\begin{center}
$a=\underset{w\in \mathbb{F}^{+}(G:a),\,u_{w}\in \{1,*\}}{\sum }%
p_{w}L_{w}^{u_{w}}$
\end{center}

\strut

in the graph $W^{*}$-probability space $\left( W^{*}(G),E\right) $.

\strut

\strut

\subsubsection{Lattice Path Model}

\strut

\strut

\strut

Throughout this section, let $G$ be a countable directed graph and let $%
\left( W^{*}(G),E\right) $ be the graph $W^{*}$-probability space over its
diagonal subalgebra $D_{G}.$ Let $w_{1},...,w_{n}\in \Bbb{F}^{+}(G)$ and let 
$L_{w_{1}}^{u_{w_{1}}}...L_{w_{n}}^{u_{w_{n}}}\in \left( W^{*}(G),E\right) $
be a $D_{G}$-valued random variable. In this section, we will define a
lattice path model for the random variable $%
L_{w_{1}}^{u_{w_{1}}}...L_{w_{n}}^{u_{w_{n}}}.$ Recall that if $%
w=e_{1}....e_{k}\in FP(G)$ with $e_{1},...,e_{k}\in E(G),$ then we can
define the length $\left| w\right| $ of $w$ by $k.$ i.e.e, the length $%
\left| w\right| $ of $w$ is the cardinality $k$ of the admissible edges $%
e_{1},...,e_{k}.$

\strut

\begin{definition}
Let $G$ be a countable directed graph and $\Bbb{F}^{+}(G),$ the free
semigroupoid. If $w\in \Bbb{F}^{+}(G),$ then $L_{w}$ is the corresponding $%
D_{G}$-valued random variable in $\left( W^{*}(G),E\right) .$ We define the
lattice path $l_{w}$ of $L_{w}$ and the lattice path $l_{w}^{-1}$ of $%
L_{w}^{*}$ by the lattice paths satisfying that ;

\strut 

(i) \ \ the lattice path $l_{w}$ starts from $*=(0,0)$ on the $\Bbb{R}^{2}$%
-plane.

\strut 

(ii) \ if $w\in V(G),$ then $l_{w}$ has its end point $(0,1).$

\strut 

(iii) if $w\in E(G),$ then $l_{w}$ has its end point $(1,1).$

\strut 

(iv) if $w\in E(G),$ then $l_{w}^{-1}$ has its end point $(-1,-1).$

\strut 

(v) \ if $w\in FP(G)$ with $\left| w\right| =k,$ then $l_{w}$ has its end
point $(k,k).$

\strut 

(vi) if $w\in FP(G)$ with $\left| w\right| =k,$ then $l_{w}^{-1}$ has its
end point $(-k,-k).$

\strut 

Assume that finite paths $w_{1},...,w_{s}$ in $FP(G)$ satisfy that $%
w_{1}...w_{s}\in FP(G).$ Define the lattice path $l_{w_{1}...w_{s}}$ by the
connected lattice path of the lattice paths $l_{w_{1}},$ ..., $l_{w_{s}}.$
i.e.e, $l_{w_{2}}$ starts from $(k_{w_{1}},k_{w_{1}})\in \Bbb{R}^{+}$ and
ends at $(k_{w_{1}}+k_{w_{2}},k_{w_{1}}+k_{w_{2}}),$ where $\left|
w_{1}\right| =k_{w_{1}}$ and $\left| w_{2}\right| =k_{w_{2}}.$ Similarly, we
can define the lattice path $l_{w_{1}...w_{s}}^{-1}$ as the connected path
of $l_{w_{s}}^{-1},$ $l_{w_{s-1}}^{-1},$ ..., $l_{w_{1}}^{-1}.$
\end{definition}

\strut

\begin{definition}
Let $G$ be a countable directed graph and assume that $%
L_{w_{1}},...,L_{w_{n}}$ are generators of $\left( W^{*}(G),E\right) .$ Then
we have the lattice paths $l_{w_{1}},$ ..., $l_{w_{n}}$ of $L_{w_{1}},$ ..., 
$L_{w_{n}},$ respectively in $\Bbb{R}^{2}.$ Suppose that $%
L_{w_{1}}^{u_{w_{1}}}...L_{w_{n}}^{u_{w_{n}}}\neq 0_{D_{G}}$ in $\left(
W^{*}(G),E\right) ,$ where $u_{w_{1}},...,u_{w_{n}}\in \{1,*\}.$ Define the
lattice path $l_{w_{1},...,w_{n}}^{u_{w_{1}},...,u_{w_{n}}}$ of nonzero $%
L_{w_{1}}^{u_{w_{1}}}...L_{w_{n}}^{u_{w_{n}}}$ by the connected lattice path
of $l_{w_{1}}^{t_{w_{1}}},$ ..., $l_{w_{n}}^{t_{w_{n}}},$ where $t_{w_{j}}=1$
if $u_{w_{j}}=1$ and $t_{w_{j}}=-1$ if $u_{w_{j}}=*.$ Assume that $%
L_{w_{1}}^{u_{w_{1}}}...L_{w_{n}}^{u_{w_{n}}}$ $=$ $0_{D_{G}}.$ Then the
empty set $\emptyset $ in $\Bbb{R}^{2}$ is the lattice path of it. We call
it the empty lattice path. By $LP_{n},$ we will denote the set of all
lattice paths of the $D_{G}$-valued random variables having their forms of $%
L_{w_{1}}^{u_{w_{1}}}...L_{w_{n}}^{u_{w_{n}}},$ including empty lattice path.
\end{definition}

\strut

Also, we will define the following important property on the set of all
lattice paths ;

\strut

\begin{definition}
Let $l_{w_{1},...,w_{n}}^{u_{w_{1}},...,u_{w_{n}}}\neq \emptyset $ be a
lattice path of $L_{w_{1}}^{u_{w_{1}}}...L_{w_{n}}^{u_{w_{n}}}\neq 0_{D_{G}}$
in $LP_{n}.$ If the lattice path $%
l_{w_{1},...,w_{n}}^{u_{w_{1}},...,u_{w_{n}}}$ starts from $*$ and ends on
the $*$-axis in $\Bbb{R}^{+},$ then we say that the lattice path $%
l_{w_{1},...,w_{n}}^{u_{w_{1}},...,u_{w_{n}}}$ has the $*$-axis-property. By 
$LP_{n}^{*},$ we will denote the set of all lattice paths having their forms
of $l_{w_{1},...,w_{n}}^{u_{w_{1}},...,u_{w_{n}}}$ which have the $*$%
-axis-property. By little abuse of notation, sometimes, we will say that the 
$D_{G}$-valued random variable $L_{w_{1}}^{u_{w_{1}}}...L_{w_{n}}^{u_{w_{n}}}
$satisfies the $*$-axis-property if the lattice path $%
l_{w_{1},...,w_{n}}^{u_{w_{1}},...,u_{w_{n}}}$ of it has the $*$%
-axis-property.
\end{definition}

\strut

The following theorem shows that finding $E\left(
L_{w_{1}}^{u_{w_{1}}}...L_{w_{n}}^{u_{w_{n}}}\right) $ is checking the $*$%
-axis-property of $L_{w_{1}}^{u_{w_{1}}}...L_{w_{n}}^{u_{w_{n}}}.$

\strut

\begin{theorem}
(See [15]) Let $L_{w_{1}}^{u_{w_{1}}}...L_{w_{n}}^{u_{w_{n}}}\in \left(
W^{*}(G),E\right) $ be a $D_{G}$-valued random variable, where $%
u_{w_{1}},...,u_{w_{n}}\in \{1,*\}.$ Then $E\left(
L_{w_{1}}^{u_{w_{1}}}...L_{w_{n}}^{u_{w_{n}}}\right) $ $\neq $ $0_{D_{G}}$
if and only if $L_{w_{1}}^{u_{w_{1}}}...L_{w_{n}}^{u_{w_{n}}}$ has the $*$%
-axis-property (i.e., the corresponding lattice path $%
l_{w_{1},...,w_{n}}^{u_{w_{1}},...,u_{w_{n}}}$ of $%
L_{w_{1}}^{u_{w_{1}}}...L_{w_{n}}^{u_{w_{n}}}$ is contained in $LP_{n}^{*}.$
Notice that $\emptyset \notin LP_{n}^{*}.$) \ $\square $
\end{theorem}

\strut \strut

By the previous theorem, we can conclude that $E\left(
L_{w_{1}}^{u_{w_{1}}}...L_{w_{n}}^{u_{w_{n}}}\right) =L_{v},$ for some $v\in
V(G)$ if and only if the lattice path $%
l_{w_{1},...,w_{n}}^{u_{w_{1}},...,u_{w_{n}}}$ has the $*$-axis-property
(i.e., $l_{w_{1},...,w_{n}}^{u_{w_{1}},...,u_{w_{n}}}\in LP_{n}^{*}$).\strut
\strut \strut

\strut

\strut \strut

\subsubsection{$D_{G}$-Valued Moments and Cumulants of Random Variables\strut
}

\bigskip

\bigskip

Let $w_{1},...,w_{n}\in \mathbb{F}^{+}(G)$, $u_{1},...,u_{n}\in \{1,*\}$ and
let $L_{w_{1}}^{u_{1}}...L_{w_{n}}^{u_{n}}\in \left( W^{*}(G),E\right) $ be
a $D_{G}$-valued random variable. Recall that, in the previous section, we
observed that the $D_{G}$-valued random variable $%
L_{w_{1}}^{u_{1}}...L_{w_{n}}^{u_{n}}=L_{v}\in \left( W^{*}(G),E\right) $
with $v\in V(G)$ if and only if the lattice path $%
l_{w_{1},...,w_{n}}^{u_{1},...,u_{n}}$ of \ $%
L_{w_{1}}^{u_{1}}...L_{w_{n}}^{u_{n}}$ has the $*$-axis-property
(equivalently, $l_{w_{1},...,w_{n}}^{u_{1},...,u_{n}}\in LP_{n}^{*}$).
Throughout this section, fix a $D_{G}$-valued random variable $a\in \left(
W^{*}(G),E\right) .$ Then the $D_{G}$-valued random variable $a$ has the
following Fourier expansion,

\bigskip

\begin{center}
$a=\underset{v\in V(G:a)}{\sum }p_{v}L_{v}+\underset{l\in FP_{*}(G:a)}{\sum }%
\left( p_{l}L_{l}+p_{l^{t}}L_{l}\right) +\underset{w\in
FP_{*}^{c}(G:a),~u_{w}\in \{1,*\}}{\sum }p_{w}L_{w}^{u_{w}}.$
\end{center}

\bigskip

Let's observe the new $D_{G}$-valued random variable $d_{1}ad_{2}a...d_{n}a%
\in \left( W^{*}(G),E\right) ,$ where $d_{1},...,d_{n}\in D_{G}$ and $a\in
W^{*}(G)$ is given. Put

\strut

\begin{center}
$d_{j}=\underset{v_{j}\in V(G:d_{j})}{\sum }q_{v_{j}}L_{v_{j}}\in D_{G},$
for \ $j=1,...,n.$
\end{center}

\strut

Notice that $V(G:d_{j})=\mathbb{F}^{+}(G:d_{j}),$ since $d_{j}\in
D_{G}\hookrightarrow W^{\ast }(G).$ Then

\strut \strut

\begin{proposition}
(See [16]) Let $a\in \left( W^{*}(G),E\right) $ be given as above. Then the $%
n$-th moment of $a$ is

\strut 

$\ \ E\left( d_{1}a...d_{n}a\right) =\underset{(v_{1},...,v_{n})\in \Pi
_{j=1}^{n}V(G:d_{j})}{\sum }\left( \Pi _{j=1}^{n}q_{v_{j}}\right) $

\strut 

$\ \ \ \underset{(w_{1},...,w_{n})\in \mathbb{F}^{+}(G:a)^{n},\,u_{w_{j}}\in
\{1,*\},\,l_{w_{1},...,w_{n}}^{u_{w_{1}},...,u_{w_{n}}}\in LP_{n}^{*}}{\sum }%
\left( \Pi _{j=1}^{n}p_{w_{j}}\right) $

\strut 

$\ \ \ \ \ \ \ \ \ \ \ \ \ \ \ \left( \Pi _{j=1}^{n}\delta
_{(v_{j},x_{j},y_{j}:u_{w_{j}})}\right) \,\,E\left(
L_{w_{1}}^{u_{w_{1}}}...L_{w_{n}}^{u_{w_{n}}}\right) ,$

where

$\ \ \ \ \ \ \ \ \ \ \ \delta _{(v_{j},x_{j},y_{j})}=\left\{ 
\begin{array}{lll}
\delta _{v_{j},x_{j}} &  & \text{if }u_{j}=1 \\ 
\delta _{v_{j},y_{j}} &  & \text{if }u_{j}=*
\end{array}
\right. $

$\square $
\end{proposition}

\strut

Let $w_{1},...,w_{n}\in FP(G)$ be finite paths and \ $u_{1},...,u_{n}\in
\{1,*\}$. Then, by the M\"{o}bius inversion, we have

\strut

(1.13)$\ \ $

\begin{center}
$k_{n}\left( L_{w_{1}}^{u_{1}}~,...,~L_{w_{n}}^{u_{n}}\right) =\underset{\pi
\in NC(n)}{\sum }\widehat{E}(\pi )\left( L_{w_{1}}^{u_{1}}~\otimes
...\otimes ~L_{w_{n}}^{u_{n}}\right) \mu (\pi ,1_{n}),$
\end{center}

\strut

where $\widehat{E}=\left( E^{(n)}\right) _{n=1}^{\infty }$ is the moment
multiplicative bimodule map induced by the conditional expectation $E$ (See
[16]) and where $NC(n)$ is the collection of all noncrossing partition over $%
\{1,...,n\}.$ 

\strut

\begin{definition}
Let $NC(n)$ be the set of all \ noncrossing partition over $\{1,...,n\}$ and
let $L_{w_{1}}^{u_{1}},$ $...,$ $L_{w_{n}}^{u_{n}}\in \left(
W^{*}(G),E\right) $ be $D_{G}$-valued random variables, where $%
u_{1},...,u_{n}\in \{1,*\}.$ We say that the $D_{G}$-valued random variable $%
L_{w_{1}}^{u_{1}}...L_{w_{n}}^{u_{n}}$ is $\pi $-connected if the $\pi $%
-dependent $D_{G}$-moment of it is nonvanishing, for $\pi \in NC(n).$ In
other words, the random variable $L_{w_{1}}^{u_{1}}...L_{w_{n}}^{u_{n}}$ is $%
\pi $-connected, for $\pi \in NC(n),$ if

\strut 

$\ \ \ \ \ \ \ \ \ \ \ \widehat{E}(\pi )\left( L_{w_{1}}^{u_{1}}~\otimes
...\otimes ~L_{w_{n}}^{u_{n}}\right) \neq 0_{D_{G}}.$

\strut 

i.e., there exists a vertex $v\in V(G)$ such that

\strut 

$\ \ \ \ \ \ \ \ \ \ \ \widehat{E}(\pi )\left( L_{w_{1}}^{u_{1}}~\otimes
...\otimes ~L_{w_{n}}^{u_{n}}\right) =L_{v}.$
\end{definition}

\strut

For convenience, we will define the following subset of $NC(n)$ ;

\strut

\begin{definition}
Let $NC(n)$ be the set of all noncrossing partitions over $\{1,...,n\}$ and
fix a $D_{G}$-valued random variable $L_{w_{1}}^{u_{1}}...L_{w_{n}}^{u_{n}}$
in $\left( W^{*}(G),E\right) ,$ where $u_{1},$ ..., $u_{n}\in \{1,*\}.$ For
the fixed $D_{G}$-valued random variable $%
L_{w_{1}}^{u_{1}}...L_{w_{n}}^{u_{n}},$define

\strut 

$\ \ \ C_{w_{1},...,w_{n}}^{u_{1},...,u_{n}}\overset{def}{=}\{\pi \in
NC(n):L_{w_{1}}^{u_{1}}...L_{w_{n}}^{u_{n}}$ is $\pi $-connected$\},$

\strut 

in $NC(n).$ Let $\mu $ be the M\"{o}bius function in the incidence algebra $%
I_{2}.$ Define the number $\mu _{w_{1},...,w_{n}}^{u_{1},...,u_{n}},$ for
the fixed $D_{G}$-valued random variable $%
L_{w_{1}}^{u_{1}}...L_{w_{n}}^{u_{n}},$ by

\strut 

$\ \ \ \ \ \ \ \ \ \ \ \ \ \ \ \mu _{w_{1},...,w_{n}}^{u_{1},...,u_{n}}%
\overset{def}{=}\underset{\pi \in C_{w_{1},...,w_{n}}^{u_{1},...,u_{n}}}{%
\sum }\mu (\pi ,1_{n}).$
\end{definition}

\bigskip \strut

Assume that there exists $\pi \in NC(n)$ such that $%
L_{w_{1}}^{u_{1}}...L_{w_{n}}^{u_{n}}=L_{v}$ is $\pi $-connected. Then $\pi
\in C_{w_{1},...,w_{n}}^{u_{1},...,u_{n}}$ and there exists the maximal
partition $\pi _{0}\in C_{w_{1},...,w_{n}}^{u_{1},...,u_{n}}$ such that $%
L_{w_{1}}^{u_{1}}...L_{w_{n}}^{u_{n}}=L_{v}$ is $\pi _{0}$-connected. Notice
that $1_{n}\in C_{w_{1},...,w_{n}}^{u_{1},...,u_{n}}.$ Therefore, the
maximal partition in $C_{w_{1},...,w_{n}}^{u_{1},...,u_{n}}$ is $1_{n}.$
Hence we have that ;

\bigskip \strut

\begin{lemma}
(See [16]) Let $L_{w_{1}}^{u_{1}}...L_{w_{n}}^{u_{n}}\in \left(
W^{*}(G),E\right) $ be a $D_{G}$-valued random variable having the $*$%
-axis-property. Then

\strut 

$\ \ \ \ \ \ \ \ \ \ \ E\left( L_{w_{1}}^{u_{1}}...L_{w_{n}}^{u_{n}}~\right)
=\ \widehat{E}(\pi )\left( L_{w_{1}}^{u_{1}}~\otimes ...\otimes
~L_{w_{n}}^{u_{n}}\right) ,$

\strut 

for all $\pi \in C_{w_{1},...,w_{n}}^{u_{1},...,u_{n}}.$ \ $\square $
\end{lemma}

\bigskip \strut \strut \strut

By the previous lemmas, we have that

\strut

\begin{theorem}
(See [16]) Let $n\in 2\mathbb{N}$ and let $%
L_{w_{1}}^{u_{1}},...,L_{w_{n}}^{u_{n}}\in \left( W^{*}(G),E\right) $ be $%
D_{G}$-valued random variables, where $w_{1},...,w_{n}\in FP(G)$ and $%
u_{j}\in \{1,*\},$ $j=1,...,n.$ Then

\strut 

$\ \ \ \ \ \ \ k_{n}\left( L_{w_{1}}^{u_{1}}...L_{w_{n}}^{u_{n}}~\right)
=\mu _{w_{1},...,w_{n}}^{u_{1},...,u_{n}}\cdot
E(L_{w_{1}}^{u_{1}},...,L_{w_{n}}^{u_{n}}),$

\strut 

where $\mu _{w_{1},...,w_{n}}^{u_{1},...,u_{n}}=\underset{\pi \in
C_{w_{1},...,w_{n}}^{u_{1},...,u_{n}}}{\sum }\mu (\pi ,1_{n}).$ \ $\square $
\end{theorem}

\bigskip \strut \strut \strut \strut \strut

\strut \strut

\strut \strut

\subsection{\strut $D_{G}$-Freeness on $\left( W^{*}(G),E\right) $}

\strut

\strut

Now, we will introduce the diagram-distinctness of finite paths ;

\strut

\begin{definition}
(\textbf{Diagram-Distinctness}) We will say that the finite paths $w_{1}$
and $w_{2}$ are \textbf{diagram-distinct} if $w_{1}$ and $w_{2}$ have
different diagrams in the graph $G.$ Let $X_{1}$ and $X_{2}$ be subsets of $%
FP(G).$ The subsets $X_{1}$ and $X_{2}$ are said to be diagram-distinct if $%
x_{1}$ and $x_{2}$ are diagram-distinct for all pairs $(x_{1},x_{2})$ $\in $ 
$X_{1}\times X_{2}.$
\end{definition}

\strut \strut

In [16], we found the $D_{G}$-freeness characterization on the generator set
of $W^{*}(G),$ as follows ;

\strut \strut \strut

\begin{theorem}
(See [16]) Let $w_{1},w_{2}\in FP(G)$ be finite paths. The $D_{G}$-valued
random variables $L_{w_{1}}$ and $L_{w_{2}}$ in $\left( W^{*}(G),E\right) $
are free over $D_{G}$ if and only if $w_{1}$ and $w_{2}$ are
diagram-distinct. $\square $
\end{theorem}

\strut \strut \strut

Let $a$ and $b$ be the given $D_{G}$-valued random variables. We can get the
necessary condition for the $D_{G}$-freeness of $a$ and $b,$ in terms of
their supports. \strut \strut \strut Recall that we say that the two subsets 
$X_{1}$ and $X_{2}$ of $FP(G)$ are said to be diagram-distinct if $x_{1}$
and $x_{2}$ are diagram-distinct, for all pairs $(x_{1},x_{2})$ $\in $ $%
X_{1} $ $\times $ $X_{2}.$

\strut \strut

\begin{proposition}
(See [16]) Let $a,b\in \left( W^{*}(G),E\right) $ be $D_{G}$-valued random
variables with their supports $\mathbb{F}^{+}(G:a)$ and $\mathbb{F}^{+}(G:b).
$ The $D_{G}$-valued random variables $a$ and $b$ are free over $D_{G}$ in $%
\left( W^{*}(G),E\right) $ if $FP(G:a_{1})$ and $FP(G:a_{2})$ are
diagram-distinct. $\square $
\end{proposition}

\strut \strut 

\strut \strut \strut 

\strut \strut 

\section{Diagonal Compressed Random Variables in $\left( W^{*}(G),E\right) $}

\strut

\strut

Let $G$ be a countable directed graph and $\left( W^{*}(G),E\right) ,$ the
graph $W^{*}$-probability space over the diagonal subalgebra $D_{G}.$ In
[18], we observed the diagonal compressed $D_{G}$-valued free probability on 
$\left( W^{*}(G),E\right) .$ Fix a finite subset $V=\{v_{1},...,v_{N}\}$ of
the vertex set $V(G)$ and define the diagonal compression $%
P_{V}:W^{*}(G)\rightarrow W^{*}(G)$ by

\strut

\begin{center}
$P_{V}(a)=L_{v_{1}}aL_{v_{1}}+...+L_{v_{N}}aL_{v_{N}}\in W^{*}(G),$ $\forall
a\in \left( W^{*}(G),E\right) .$
\end{center}

\strut

If the given subset $V$ is a singleton set, then we will call this diagonal
compression the vertex-compression. Let $a\in \left( W^{*}(G),E\right) $ be
a $D_{G}$-valued random variable having its expression

\strut

\begin{center}
$a=\underset{v\in V(G:a)}{\sum }p_{v}L_{v}+\underset{l\in FP_{*}(G:a)}{\sum }%
\left( p_{l}L_{l}+p_{l^{t}}L_{l}^{*}\right) +\underset{w\in
FP_{*}^{c}(G:a),\,u_{w}\in \{1,*\}}{\sum }p_{w}L_{w}^{u_{w}}.$
\end{center}

\strut

Then we have the following Fourier-like expression of the diagonal
compressed random variable $P_{V}(a)$ of $a$ by $V$ ;

\strut

$\ \ \ \ \ \ \ P_{V}(a)=\underset{v\in V(G:P_{V}(a))}{\sum }p_{v}L_{v}+%
\underset{w\in FP_{*}(G:P_{V}(a))}{\sum }\left(
p_{w}L_{w}+p_{w^{t}}L_{w}^{*}\right) $

$\ \ \ \ \ \ \ \ \ \ \ \ \ \ \ \ \ \ \ \ \ \ \ \ \ \ \ \ \ \ \ \ \ \ \ \ \ \
\ \ \ \ \ \ \ +\underset{w\in FP_{*}^{c}(G:P_{V}(a)),\,u_{w}\in \{1,*\}}{%
\sum }p_{w}L_{w}^{u_{w}},$

\strut

with

\strut

$\ \ \ \ \ \ \ \ V(G:P_{V}(a))=V\cap V(G:a),$

\strut

$\ \ \ \ \ \ \ \ FP_{*}(G:P_{V}(a))=\left( \cup
_{j=1}^{N}loop_{v_{j}}(G:a)\right) \cap FP_{*}(G:a)$

and

$\ \ \ \ \ \ \ \ FP_{*}^{c}(G:P_{V}(a))=\left( \cup
_{j=1}^{N}loop_{v_{j}}(G:a)\right) \cap FP_{*}^{c}(G:a),$

\strut

where\strut

\begin{center}
$loop_{v_{j}}(G:a)=\{l\in loop(G:a):l=v_{j}lv_{j}\},$
\end{center}

for \ $j=1,...,N.$

\strut

Hence we can apply all amalgamated free probability information on this
diagonal compressed case. In particular, we can compute the $D_{G}$-valued
moments and $D_{G}$-valued cumulants of the diagonal compressed random
variables like Section 1.3 and hence we can get $D_{G}$-valued moment series
and $D_{G}$-valued R-transforms of the diagonal compressed random variables
(See [18]). Also, by little modification of Section 1.5, we can get the
following theorems ;

\strut

\begin{theorem}
(See [18]) Let $a,b\in \left( W^{*}(G),E\right) $ be a $D_{G}$-valued random
variable and let $V=\{v_{1},...,v_{N}\}$ be a finite subset of the vertex
set $V(G).$ Let $P_{V}(a)$ and $P_{V}(b)$ be the diagonal compressed random
variable of $a$ and $b$ by $V$ in $\left( W^{*}(G),E\right) ,$ respectively.
Then

\strut

(1) If $a$ and $b$ are free over $D_{G}$ in $\left( W^{*}(G),E\right) ,$
then $P_{V}(a)$ and $P_{V}(b)$ are free over $D_{G}$ in $\left(
W^{*}(G),E\right) .$

\strut

(2) If $a$ and $b$ satisfy that

\strut

$\ \ \ \ \ \ \ \ \ \left( V\cap V(G:a)\right) \cap \left( V\cap
V(G:b)\right) =\emptyset $

\strut and

$\ \ \ \ \ \ \ \ \ \ \ loop_{v_{i}}(G:a)\cap loop_{v_{j}}(G:b)=\emptyset ,$

\strut

for all choices $(i,j)\in \{1,...,N\}^{2},$ then $P_{V}(a)$ and $P_{V}(b)$
are free over $D_{G}$ in $\left( W^{*}(G),E\right) .$ $\square $
\end{theorem}

\strut

\begin{theorem}
(See [18]) Let $a\in \left( W^{*}(G),E\right) $ be a $D_{G}$-valued random
variable and let $V_{1}=\{v_{1}^{(1)},...,v_{N_{1}}^{(1)}\}$ and $%
V_{2}=\{v_{1}^{(2)},...,v_{N_{2}}^{(2)}\}$ be finite subsets of the vertex
set $V(G).$ Suppose that

\strut

$\ \ \ \ \ \ \ \ \ \ \ loop_{v_{i}^{(1)}}(G:a)\cap
loop_{v_{j}^{(2)}}(G:a)=\emptyset ,$

\strut \strut

for all choices $(i,j)\in \{1,...,N_{1}\}\times \{1,...,N_{2}\}.$ Then the
corresponding diagonal compressed random variables $P_{V_{1}}(a)$ and $%
P_{V_{2}}(a)$ of $a$ by $V_{1}$ and $V_{2}$ are free over $D_{G}$ in $\left(
W^{*}(G),E\right) .$ \ $\square $
\end{theorem}

\strut

Therefore, we can again characterize the $D_{G}$-freeness of the diagonal
compressed random variables in $\left( W^{\ast }(G),E\right) $ by the
subsets of the free semigroupoid $\mathbb{F}^{+}(G)$ of the graph $G.$ Also,
we can get the diagonal compressed R-transform calculus like in Section 1.6.

\strut

\strut

\strut

\section{Off-Diagonal Random Variables in $\left( W^{*}(G),E\right) $}

\strut

\strut

Throughout this chapter, let $G$ be a countable directed graph and $\mathbb{F%
}^{+}(G),$ the free semigroupoid of $G$ and let $\left( W^{*}(G),E\right) $
be the graph $W^{*}$-probability space over the diagonal subalgebra $D_{G}.$
In this chapter, we will consider the off-diagonal compressed random
variables in $\left( W^{*}(G),E\right) $ over $D_{G}.$ Let's fix $v_{1}\neq
v_{2}$ in $V(G)$ and let $L_{v_{1}}$ and $L_{v_{2}}$ be the corresponding
projections in $\left( W^{*}(G),E\right) .$ Define a subset $%
FP_{v_{1}}^{v_{2}}(G)$ of $FP(G)$ by

\bigskip

\begin{center}
$FP_{v_{1}}^{v_{2}}(G)\overset{def}{=}\{w\in FP(G):w=v_{1}wv_{2}\},$
\end{center}

\bigskip

for $v_{1},v_{2}\in V(G).$

\strut

\begin{definition}
Let $v_{1}\neq v_{2}\in V(G)$ be given. For any $D_{G}$-valued random
variable $a\in \left( W^{*}(G),E\right) ,$ define the $(v_{1},v_{2})$%
-off-diagonal compressed random variable of $a$ (in short $(v_{1},v_{2})$%
-compressed random variable of $a$) $_{v_{1}}a_{v_{2}}$ by

\strut

$\ \ \ \ \ \ \ \ \ \ \ \ \ \ _{v_{1}}a_{v_{2}}\overset{def}{=}%
L_{v_{1}}aL_{v_{2}}\in \left( W^{*}(G),E\right) .$
\end{definition}

\strut

Let $a\in \left( W^{*}(G),E\right) $ be an arbitrary $D_{G}$-valued random
variable having the following Fourier expansion,

\strut

\begin{center}
$a=\underset{v\in V(G:a)}{\sum }p_{v}L_{v}+\underset{l\in FP_{*}(G:a)}{\sum }%
\left( p_{l}L_{l}+p_{l^{t}}L_{l}^{*}\right) +\underset{w\in
FP_{*}^{c}(G:a),\,u\in \{1,*\}}{\sum }p_{w}L_{w}^{u}.$
\end{center}

\strut

We will denote

\strut

\begin{center}
$a_{d}=\underset{v\in V(G:a)}{\sum }p_{v}L_{v},$ \ $a_{(*)}=\underset{l\in
FP_{*}(G:a)}{\sum }\left( p_{l}L_{l}+p_{l^{t}}L_{l}^{*}\right) $
\end{center}

and

\begin{center}
$a_{(non-*)}=\underset{w\in FP_{*}^{c}(G:a),\,u\in \{1,*\}}{\sum }%
p_{w}L_{w}^{u},$
\end{center}

\strut

for the given $D_{G}$-valued random variable $a.$ Thus the $(v_{1},v_{2})$%
-compressed random variable $_{v_{1}}a_{v_{2}}$ is determined by

\strut \strut

$\
_{v_{1}}a_{v_{2}}=L_{v_{1}}aL_{v_{2}}=L_{v_{1}}a_{d}L_{v_{2}}+L_{v_{1}}a_{(%
\ast )}L_{v_{2}}+L_{v_{1}}a_{(non-\ast )}L_{v_{2}}$

\strut

$\ \ \ \ \ \ \ \ \ =0_{D_{G}}+L_{v_{1}}a_{(\ast
)}L_{v_{2}}+L_{v_{1}}a_{(non-\ast )}L_{v_{2}}.$

\strut

By definition, we have the following partition of $FP(G:a),$ for the given
random variable $a\in $ $\left( W^{*}(G),E\right) $ ;

\strut

\begin{center}
$\{l=v_{1}lv_{2}:l\in FP_{*}(G:a)\}=FP_{*}(G:a)\cap FP_{v_{1}}^{v_{2}}(G)$
\end{center}

and

\begin{center}
$\{w=v_{1}wv_{2}:w\in FP_{*}^{c}(G:a)\}=FP_{*}^{c}(G:a)\cap
FP_{v_{1}}^{v_{2}}(G).$
\end{center}

\strut

So, we have that ;

\strut \strut

\begin{lemma}
Let $a\in \left( W^{*}(G),E\right) $ be a $D_{G}$-valued random variable and
let $v_{1}\neq v_{2}\in V(G)$ be the fixed vertices. Let $%
_{v_{2}}a_{v_{1}}\in \left( W^{*}(G),E\right) $ be the $(v_{1},v_{2})$%
-off-diagonal compressed random variable. Then

\strut

$\ \ \ _{v_{1}}a_{v_{2}}=\underset{w=v_{1}wv_{2}\in FP(G:a)}{\sum }%
p_{w}L_{w}+\underset{w^{\prime }=v_{2}w^{\prime }v_{1}\in FP(G:a)}{\sum }%
p_{w^{\prime }}L_{w^{\prime }}^{*}.$
\end{lemma}

\strut

\begin{proof}
By the relation that

\strut

$\ \ \ \ \ \ L_{w}=L_{vw}=L_{v}L_{w},$ \ \ $L_{w}=L_{wv^{\prime
}}=L_{w}L_{v^{\prime }}$

and

$\ \ \ \ \ \ L_{w}^{\ast }=L_{wv}^{\ast }=L_{v}L_{w}^{\ast },$ \ $\
L_{w}^{\ast }=L_{v^{\prime }w}^{\ast }=L_{w}^{\ast }L_{v^{\prime }},$

\strut

under the weak topology ($v,v^{\prime }\in V(G)$), we have that

\strut

$\ \ \ \ \ \ \ \ \ \ \ \ \ L_{w}=L_{vwv^{\prime }}=L_{v}L_{w}L_{v^{\prime }}$

\strut and

$\ \ \ \ \ \ \ \ \ \ \ \ \ L_{w}^{*}=L_{vwv^{\prime }}^{*}=L_{v^{\prime
}}L_{w}^{*}L_{v},$

\strut

whenever $w=vwv^{\prime }$ is a non-loop finite path, for $v,v^{\prime }\in
V(G).$ Thus, if $w\in FP(G:a),$ then we have that

\strut

$\ \ \ \ \ \ \ L_{v_{1}}L_{w}L_{v_{2}}=\left\{ 
\begin{array}{lll}
L_{w} &  & \text{if }w=v_{1}wv_{2} \\ 
&  &  \\ 
0_{D_{G}} &  & \text{otherwise}
\end{array}
\right. $

and

$\ \ \ \ \ \ \ L_{v_{1}}L_{w}^{\ast }L_{v_{2}}=\left\{ 
\begin{array}{lll}
L_{w}^{\ast } &  & \text{if }w=v_{2}wv_{1} \\ 
&  &  \\ 
0_{D_{G}} &  & \text{otherwise.}
\end{array}
\right. $

\strut

Therefore, the $(v_{1},v_{2})$-compressed random variable $_{v_{1}}a_{v_{2}}$
can be

\strut

$\ \ \ _{v_{1}}a_{v_{2}}=L_{v_{1}}aL_{v_{2}}=L_{v_{1}}\left( \underset{w\in 
\mathbb{F}^{+}(G:a),\,u_{w}\in \{1,*\}}{\sum }p_{w}L_{w}^{u_{w}}\right)
L_{v_{2}}$

\strut

$\ \ \ \ \ \ \ \ \ \ \ =L_{v_{1}}\left( \underset{w\in FP(G:a),\,u_{w}\in
\{1,\ast \}}{\sum }p_{w}L_{w}^{u_{w}}\right) L_{v_{2}}$

\strut

$\ \ \ \ \ \ \ \ \ \ \ =\underset{w\in FP(G:a),\,u_{w}\in \{1,\ast \}}{\sum }%
p_{w}\left( L_{v_{1}}L_{w}^{u_{w}}L_{v_{2}}\right) $

\strut

$\ \ \ \ \ \ \ \ \ \ \ =\underset{w=v_{1}wv_{2}\in FP(G:a)}{\sum }p_{w}L_{w}+%
\underset{w^{\prime }=v_{2}w^{\prime }v_{1}\in FP(G:a)}{\sum }p_{w^{\prime
}}L_{w^{\prime }}^{*}.$

\strut \strut \strut
\end{proof}

\strut \strut

For the convenience, we introduce the following new notation ;

\strut

\begin{quote}
\frame{\textbf{Notation}} Let $v_{1}\neq v_{2}\in V(G)$ be given as before
and let $a\in \left( W^{*}(G),E\right) $ be a $D_{G}$-valued random
variable. Define

\strut
\end{quote}

\begin{center}
$FP_{v_{1}}^{v_{2}}(G:a)_{*}\overset{def}{=}FP_{*}(G:a)\cap
FP_{v_{1}}^{v_{2}}(G)$
\end{center}

\begin{quote}
and
\end{quote}

\begin{center}
$FP_{v_{2}}^{v_{1}}(G:a)_{*}^{c}\overset{def}{=}FP_{*}^{c}(G:a)\cap
FP_{v_{2}}^{v_{1}}(G).$
\end{center}

\begin{quote}
$\square $
\end{quote}

\strut

\strut

\strut

\strut

\subsection{Off-Diagonal Compressed Moments and Cumulants}

\strut

\strut

Let $v_{1}\neq v_{2}\in V(G)$ be the fixed vertices. In this section, we
will consider the moments and cumulants of the $(v_{1},v_{2})$-off-diagonal
compressed random variables in the graph $W^{*}$-probability space $\left(
W^{*}(G),E\right) ,$ over the diagonal subalgebra $D_{G}.$ We have that

\strut

\begin{center}
$_{v_{1}}a_{v_{2}}=\underset{w=v_{1}wv_{2}\in FP_{v_{1}}^{v_{2}}(G:a)}{\sum }%
p_{w}L_{w}+\underset{w^{\prime }=v_{2}w^{\prime }v_{1}\in
FP_{v_{2}}^{v_{1}}(G:a)}{\sum }p_{w^{\prime }}L_{w^{\prime }}^{*}.$
\end{center}

\strut

So, to compute the $D_{G}$-valued moments and the $D_{G}$-valued cumulants
of $(v_{1},v_{2})$-off-diagonal compressed random variables in $\left(
W^{*}(G),E\right) $ is to compute the $D_{G}$-valued moments and the $D_{G}$%
-valued cumulants of the $D_{G}$-valued random variables $x\in \left(
W^{*}(G),E\right) $ such that

\strut

\begin{center}
$x=\underset{l_{1}=v_{1}l_{1}v_{2}\in FP(G:x)}{\sum }p_{l_{1}}L_{l_{1}}+%
\underset{l_{2}=v_{2}l_{2}v_{1}\in FP(G:x)}{\sum }p_{l_{2}}L_{l_{2}}^{*}$
\end{center}

\strut \strut

in $\left( W^{*}(G),E\right) .$

\strut

Suppose that $a\in \left( W^{*}(G),E\right) $ is a $D_{G}$-valued random
variable and assume that

\strut

\begin{center}
$FP_{*}(G:a)=\{w_{1},w_{2},...\}\subset FP(G).$
\end{center}

\strut

Then, in terms of $FP_{*}(G:a),$ the $D_{G}$-valued random variable $a$ has
the following summands

\strut

\begin{center}
$p_{w_{1}}L_{w_{1}},$ $\ \ p_{w_{1}^{t}}L_{w_{1}}^{*},$ \ $%
p_{w_{2}}L_{w_{2}},$ \ \ $p_{w_{2}^{t}}L_{w_{2}}^{*},...,$
\end{center}

\strut

where $p_{w_{j}}\,,\,\,p_{w_{j}^{t}}\in \mathbb{C}$. By the above
observation, we have the following result ;

\strut

\begin{theorem}
Let $v_{1}\neq v_{2}\in V(G)$ be the fixed vertices in the graph $G$ and let 
$x\in \left( W^{*}(G),E\right) $ be a $D_{G}$-valued random variable with
its Fourier-like expression,

\strut

$\ \ \ \ \ \ \ \ \ \underset{l_{1}=v_{1}l_{1}v_{2}\in FP(G:x)}{\sum }%
p_{l_{1}}L_{l_{1}}+\underset{l_{2}=v_{2}l_{2}v_{1}\in FP(G:x)}{\sum }%
p_{l_{2}}L_{l_{2}}^{*}.$

\strut

Then the $n$-th moments and $n$-th cumulants of $x$ vanish, for all $n\in %
\mathbb{N}.$
\end{theorem}

\strut

\begin{proof}
(1) Let $n=1.$ Then the first moments and the first cumulants of the $D_{G}$%
-valued random variable $x$ vanish ;

\strut

$\ \ \ \ \ \ \ \ \ \ \ \ \ \ \ \ \ \ \ \ \ \ \ E\left( x\right)
=k_{1}(x)=0_{D_{G}},$

\strut

since $V(G:x)=\emptyset .$

\strut

(2) Let $n>1$ in $\mathbb{N}.$ Then the $n$-th $D_{G}$-valued moments vanish
; By Section 1.4, we have that the $n$-th moment of the $D_{G}$-valued
random variable $x$ is

\strut

$E\left( d_{1}x...d_{n}x\right) $

\strut \strut

$\ =\underset{\pi \in NC(n)}{\sum }\,\underset{(v^{(1)},...,v^{(n)})\in \Pi
_{j=1}^{n}V(G:d_{j})}{\sum }\left( \Pi _{j=1}^{n}q_{v^{(j)}}\right) \,$

\strut

$\ \ \ \ \ \ \ \ \ \ \ \ \ \ \underset{(w_{1},...,w_{n})\in
FP(G:x)^{n},\,w_{j}=x_{j}w_{j}y_{j},\,u_{w_{j}}\in
\{1,*\},\,l_{w_{1},...,w_{n}}^{u_{1},...,u_{n}}\in LP_{n}^{*}}{\sum }\left(
\Pi _{j=1}^{n}p_{w_{j}}\right) $

\strut

$\ \ \ \ \ \ \ \ \ \ \ \ \ \ \ \ \ \ \,\left( \Pi _{j=1}^{n}\delta
_{(v^{(j)},\,x_{j},\,y_{j}:u_{w_{j}})}\right) E\left(
L_{w_{1}}^{u_{1}}...L_{w_{n}}^{u_{n}}\right) ,$

\strut

where $x_{j},\,y_{j}\in \{v_{1},v_{2}\}$ and where $d_{j}=\underset{%
v^{(j)}\in V(G:d_{j})}{\sum }q_{v^{(j)}}L_{v^{(j)}}\in D_{G}$ are arbitrary,
\ $j=1,...,n,$ for all $n\in \mathbb{N},$ and where $LP_{n}^{*}$ is the
lattice path model satisfying the $*$-axis-property (See Section 1.2). But
to get the nonvanishing $n$-th cumulant of $x$, we need to have at least one
summand of $d_{1}a...d_{n}a,$ $L_{w_{1}}^{u_{1}}...L_{w_{n}}^{u_{n}}=L_{v},$
for some $v\in V(G).$ Equivalently, the lattice path $%
l_{w_{1},...,w_{n}}^{u_{1},...,u_{n}}$ should have the $*$-axis-property
(i.e., $l_{w_{1},...,w_{n}}^{u_{1},...,u_{n}}\in LP_{n}^{*}$). To do that,
at least, we need to have the nonempty $FP_{*}(G:x).$ But

\strut

$\ \ \ \ \ \ \ \ \ \ \ V(G)\nsubseteq F^{+}(G:a)$ and $FP_{*}(G:x)=\emptyset
.$

\strut

\strut Therefore,

$\ \ \ \ \ \ \ \ \ \ \ \ \ \ \ \ \ k_{1}(d_{1}x)=0_{D_{G}}=E(d_{1}a)$

and

$\ \ \ \ \ \ \ \ \ \ \ k_{n}\left( a,...,a\right) =k_{n}\left(
a_{(*)},...,a_{(*)}\right) =0_{D_{G}}.$

\strut

Indeed, we have that

\strut

$\ \ \ \ \ \ \ \ \ \ \ FP(G:x)\subseteq FP_{v_{1}}^{v_{2}}(G)\cup
FP_{v_{2}}^{v_{1}}(G).$

\strut

Moreover, if $w=v_{1}wv_{2}\in FP(G:x)\cap FP_{v_{1}}^{v_{2}}(G),$ then the $%
L_{w}$-term of $x$ exists but the $L_{w}^{*}$-term does not exists in the
Fourier expansion of $x.$ Similarly, if $w^{\prime }=v_{2}wv_{1}\in
FP(G:x)\cap FP_{v_{2}}^{v_{1}}(G),$ then the $L_{w^{\prime }}^{*}$-term of $%
x $ exists but the $L_{w^{\prime }}$-term does not exists in the Fourier
expansion of $x.$ Therefore, each lattice path of $%
L_{w_{1}}^{u_{1}}...L_{w_{n}}^{u_{n}}$ does not have the $*$-axis-property.
Since all $n$-th $D_{G}$-valued cumulants of $x$ vanish, all $k$-th $D_{G}$%
-valued moments of $x$ vanish, by the M\"{o}bius inversion.
\end{proof}

\strut \strut

Consider the $D_{G}$-valued random variable

\strut

\begin{center}
$a=L_{v}+L_{w_{1}}+L_{w_{1}}^{*}+L_{w_{2}}^{*},$
\end{center}

\strut

where $v\in V(G)$ and $w_{1}\neq w_{2}\in FP(G),$ with $%
w_{1}=v_{1}w_{1}v_{2} $ and $w_{2}=v_{2}w_{2}v_{1}.$ Then the $(v_{1},v_{2})$%
-off-diagonal compressed random variable of the $D_{G}$-valued random
variable $a$ is

\strut

\begin{center}
$_{v_{1}}a_{v_{2}}=L_{w_{1}}+L_{w_{2}}^{*}.$
\end{center}

\strut

So, we have that

\strut

\begin{center}
$FP_{*}(G:\,_{v_{1}}a_{v_{2}})=\emptyset $
\end{center}

and

\begin{center}
$FP_{*}^{c}(G:\,_{v_{1}}a_{v_{2}})=\{w_{1},w_{2}\}.$
\end{center}

\strut

Therefore, the $n$-th moments and the $n$-th cumulants of $_{v_{1}}a_{v_{2}}$
vanish, for all $n\in \mathbb{N}.$

\strut

\strut

\strut

\subsection{Off-Diagonal Compressed $D_{G}$-Freeness}

\strut

\strut

\strut

Let $G$ be a countable directed graph and $\mathbb{F}^{+}(G),$ the free
semigroupoid of $G$ and let $\left( W^{*}(G),E\right) $ be the corresponding
graph $W^{*}$-probability space over the diagonal subalgebra $D_{G}.$ In
this section, we will consider the $D_{G}$-freeness of the off-diagonal
compressed random variables. Throughout this section, let $v_{1},v_{2},v_{3}$
and $v_{4}$ be mutually distinct vertices in $V(G).$ Since $%
FP_{*}(G:\,_{v_{1}}a_{v_{2}})=\emptyset ,$ it has vanishing $D_{G}$-valued $%
n $-th moments and $n$-th cumulants, for all $n\in \mathbb{N}.$ Therefore,
automatically, the $D_{G}$-valued moments vanish, by [12].)

\strut

\begin{proposition}
Let $a$ and $b$ be $D_{G}$-valued random variables in the graph $W^{*}$%
-probability space $\left( W^{*}(G),E\right) $ and let $v_{1}\neq v_{2}$ be
the given vertices in $V(G).$ If $a$ and $b$ have the diagram-distinct
supports, then the $(v_{1},v_{2})$-off-diagonal compressed random variables $%
_{v_{1}}a_{v_{2}}\equiv L_{v_{1}}aL_{v_{2}}$ and $_{v_{1}}b_{v_{2}}\equiv
L_{v_{1}}bL_{v_{2}}$ are free over $D_{G}$ in $\left( W^{*}(G),E\right) .$
\end{proposition}

\strut

\begin{proof}
By the diagram-distinctness of $FP(G:a)$ and $FP(G:b),$ $loop^{c}(G:a)$ and $%
loop^{c}(G:b)$ are diagram-distinct, too. Note that

\strut

$\ \ \ \ \ \ \ \ \ FP\left( G:\,_{v_{1}}a_{v_{2}}\right) \subset
loop^{c}(G:a)$

and

$\ \ \ \ \ \ \ \ \ FP\left( G:\,_{v_{1}}b_{v_{2}}\right) \subset
loop^{c}(G:b).$

\strut

Therefore, $_{v_{1}}a_{v_{2}}$ and $_{v_{1}}b_{v_{2}}$ are free over $D_{G}.$
\end{proof}

\strut

\strut Now, we will consider the other case ;

\strut

\begin{proposition}
Let $v_{1},v_{2},v_{3}$ and $v_{4}$ be the mutually distinct vertices in $%
V(G)$ and let $a\in \left( W^{*}(G),E\right) $ be a $D_{G}$-valued random
variable. Then the $(v_{1},v_{2})$-off-diagonal compressed random variable $%
_{v_{1}}a_{v_{2}}$ and the $(v_{3},v_{4})$-off-diagonal compressed random
variable $_{v_{3}}a_{v_{4}}$ are free over $D_{G}$ in $\left(
W^{*}(G),E\right) .$
\end{proposition}

\strut

\begin{proof}
Since $v_{1},...,v_{4}$ are mutually distinct, the $(v_{1},v_{2})$-off
diagonal compressed random variable $_{v_{1}}a_{v_{2}}$ and the $%
(v_{3},v_{4})$-off-diagonal compressed random variable $_{v_{3}}a_{v_{4}}$
have the diagram-distinct supports. Therefore, they are free over $D_{G}$ in 
$\left( W^{*}(G),E\right) .$
\end{proof}

\strut

\begin{corollary}
Let $v_{1},v_{2},v_{3}$ and $v_{4}$ be the given vertices in $V(G).$ Define
two subsets

\strut

$\ \ \ \ \ \ W^{*}(G)_{v_{1}}^{v_{2}}=L_{v_{1}}W^{*}(G)L_{v_{2}}$ \ and \ $%
W^{*}(G)_{v_{3}}^{v_{4}}=L_{v_{3}}W^{*}(G)L_{v_{4}},$

\strut

in the graph $W^{*}$-probability space $\left( W^{*}(G),E\right) .$ Then
these subsets are free over $D_{G}$ in $\left( W^{*}(G),E\right) .$ $\square 
$
\end{corollary}

\strut

We can regard $W^{\ast }(G)_{v_{i}}^{v_{j}}$ as an off-diagonal block of $%
W^{\ast }(G).$

\strut

\strut \strut

\strut \strut

\section{Compressed Free Probability on $\left( W^{*}(G),E\right) $}

\strut

\strut

Throughout this chapter, let $G$ be a countable directed graph and $\mathbb{F%
}^{+}(G),$ the free semigroupoid of the graph $G$ and let $\left(
W^{*}(G),E\right) $ be the graph $W^{*}$-probability space over the diagonal
subalgebra $D_{G}.$ In this chapter, we will consider the compressed random
variable $PaP$ of the $D_{G}$-valued random variable $a\in \left(
W^{*}(G),E\right) $ by the projection $P\in W^{*}(G).$ Let $v_{1},...,v_{N}$
be vertices in $V(G)$ and define the projection

\strut

\begin{center}
$P=L_{v_{1}}+...+L_{v_{N}}\in W^{*}(G).$
\end{center}

\strut

Then it is indeed a projection in $W^{*}(G)$. From now, fix the finite
vertices $v_{1},...,v_{N}\in V(G)$ and the corresponding projection $P$ $=$ $%
L_{v_{1}}$ $+$ $...$ $+$ $L_{v_{N}}$ in $W^{*}(G).$ Let $a\in \left(
W^{*}(G),E\right) $ be a $D_{G}$-valued random variable. Then naturally, we
can construct the compressed random variable $PaP$ of $a$ by $P.$ Then this
compressed random variable is again a $D_{G}$-valued random variable in $%
\left( W^{*}(G),E\right) .$

\strut

Notice that the diagonal-compressed random variable $P_{V}(a)$ by the
diagonal compression $P_{V}:W^{*}(G)\rightarrow W^{*}(G),$ for the fixed
vertex-subset $V=\{v_{1},...,v_{N}\}$ (See Chapter 2 and [18]) and the
compressed random variable $PaP$ are totally different in $\left(
W^{*}(G),E\right) .$ For example, if $a$ is a $D_{G}$-valued random variable
in $\left( W^{*}(G),E\right) ,$ then the diagonal compressed random variable
is

\strut

\begin{center}
$P_{V}(a)\overset{def}{=}L_{v_{1}}aL_{v_{1}}+...+L_{v_{N}}aL_{v_{N}}$
\end{center}

\strut

but the compressed random variable by the projection $P$ is

\strut

\begin{center}
$PaP=\left( \sum_{i=1}^{N}L_{v_{j}}\right) \,a\,\left(
\sum_{j=1}^{N}L_{v_{j}}\right) =\underset{(i,\,j)\in \{1,...,N\}^{2}}{\sum }%
L_{v_{i}}\,a\,L_{v_{j}}.$
\end{center}

\strut

Therefore, we can say that the compressed random variable $PaP$ of $a$ by $P$
satisfies that

\strut

\begin{center}
$PaP=P_{V}(a)+\underset{(i,\,j)\in \{1,...,N\}^{2},\,\,i\neq j}{\sum }%
L_{v_{i}}aL_{v_{j}}.$
\end{center}

\strut

However, in this chapter, we will observe that $P_{V}(a)$ and $PaP$ have the
same free probability information.

\strut

Again, remark that the compressed random variable $PaP$ is the sum of
diagonal compressed random variable $P_{V}(a)$ and the (sum of $D_{G}$-free)
off-diagonal compressed random variables $L_{v_{i}}aL_{v_{j}}$ \ ($i\neq j$
in $\{1,...,N\}$). Define

\strut

\begin{center}
$P_{V}^{c}(a)=\underset{(i,j)\in \{1,...,N\}^{2},\text{ }i\neq j}{\sum }%
L_{v_{i}}aL_{v_{j}}.$
\end{center}

\strut \strut \strut \strut

\begin{proposition}
Let $a\in \left( W^{*}(G),E\right) $ be a $D_{G}$-valued random variable and
let $P=\sum_{j=1}^{N}L_{v_{j}}\in W^{*}(G)$ be a projection, where $V=$ $%
\{v_{1},...,v_{N}\}$ $\subset $ $V(G)$ is the finite subset of $V(G).$ Then
the compressed random variable $PaP$ of $a$ by $P$ is

\strut

$\ \ \ \ \ \ \ \ \ \ \ \ \ \ \ PaP=P_{V}(a)+P_{V}^{c}(a),$

\strut

where $P_{V}(a)$ is the diagonal compressed random variable of $a$ by $V$ and

\strut

$\ \ \ \ \ \ \ \ \ \ \ P_{V}^{c}(a)=\underset{(i,j)\in \{1,...,N\}^{2},\text{
}i\neq j}{\sum }L_{v_{i}}aL_{v_{j}},$

\strut

In particular, $P_{V}(a)$ and $P_{V}^{c}(a)$ are free over $D_{G}$ in $%
\left( W^{*}(G),E\right) .$
\end{proposition}

\strut

\begin{proof}
By the previous discussion, the compressed random variable of $a$ by the
projection $P$ satisfies that

\strut

$\ \ \ \ \ \ \ \ \ \ \ \ \ \ \ \ \ \ \ PaP=P_{V}(a)+P_{V}^{c}(a),$

\strut where

$\ \ \ \ \ \ \ \ \ \ \ \ \
P_{V}(a)=L_{v_{1}}aL_{v_{1}}+...+L_{v_{N}}aL_{v_{N}}$

and

$\ \ \ \ \ \ \ \ \ \ \ \ \ P_{V}^{c}(a)=\underset{(i,j)\in \{1,...,N\}^{2},%
\text{ }i\neq j}{\sum }L_{v_{i}}aL_{v_{j}}.$

\strut

Then, by [18] and by Chapter 3,

\strut

$\ \ \ \ \ \ \ \ \ \ \ \Bbb{F}^{+}(G:P_{V}(a))\subset \left( V\cup
loop(G)\right) $

and

$\ \ \ \ \ \ \ \ \ \ \ \ \ \ \ \Bbb{F}^{+}(G:P_{V}^{c}(a))\subset
loop^{c}(G).$

\strut

Therefore, the supports of $P_{V}(a)$ and $P_{V}^{c}(a)$ are
diagram-distinct and hence they are free over $D_{G}.$
\end{proof}

\strut \strut

\strut

\strut \strut

\subsection{Amalgamated Moments and Cumulants of Compressed Random Variables}

\strut

\strut

Remark that the compressed random variable $PaP$ of a $D_{G}$-valued random
variable $a$ $\in $ $\left( W^{*}(G),E\right) $ by the projection $P=$ $%
\sum_{j=1}^{N}L_{v_{j}}$ $\in $ $W^{*}(G)$ has the form of

\strut

\begin{center}
$PaP=P_{V}(a)+P_{V}^{c}(a),$
\end{center}

\strut

where $V=\{v_{1},...,v_{N}\}\subset V(G)$ is the finite subset. Futhermore,
by the previous proposition, as $D_{G}$-valued random variables in $\left(
W^{\ast }(G),E\right) ,$ the diagonal compressed part $P_{V}(a)$ of $a$ and
the off-diagonal compressed part $P_{V}^{c}(a)$ are free over $D_{G}$ in $%
\left( W^{\ast }(G),E\right) .$ Therefore, we can get the following result ;

\strut

\begin{theorem}
Let $V=\{v_{1},...,v_{N}\}$ be the finite subset of the vertex set $V(G)$
and let $P=\sum_{j=1}^{N}L_{v_{j}}\in W^{*}(G)$ be the corresponding
projection. Let $a\in \left( W^{*}(G),E\right) $ be a $D_{G}$-valued random
variable and $PaP,$ the compressed random variable of $a$ by $P.$ Then the $%
n $-th cumulants of $PaP$ is

\strut

$\ \ \ \ \ \ \ \ \ \ \ k_{1}\left( d_{1}PaP\right) =\underset{v\in V\cap
(V(G:d_{1})\cap V(G:a))}{\sum }\left( q_{v}p_{v}\right) L_{v}$

\strut and

\strut \strut

$k_{n}\left( \underset{n-times}{\underbrace{d_{1}PaP,....,d_{n}PaP}}\right) $

\strut

$\ \ \ \ \ =\underset{(v^{(1)},...,v^{(n)})\in \Pi _{k=1}^{n}V(G:d_{k})}{%
\sum }\,\left( \Pi _{j=1}^{n}q_{v^{(j)}}\right) $

\strut \strut

$\ \ \ \ \ \ \ \underset{(w_{1},...,w_{n})\in \left( \left( \cup
_{k=1}^{n}loop_{v_{k}}(G:a)\right) \cup V\cap V(G:a)3\right)
^{n},~w_{j}=x_{j}w_{j}x_{j},\,%
\,l_{w_{1},...,w_{n}}^{u_{w_{1}},...,u_{w_{n}}}\in LP_{n}^{*}}{\sum }$

\strut

$\ \ \ \ \ \ \ \ \ \ \left( \Pi _{j=1}^{n}p_{w_{j}}\right) \ \ \left( \Pi
_{k=1}^{n}\delta _{v^{(k)},~x_{k}}\right) \,\,\,\ \mu
_{w_{1},...,w_{n}}^{u_{w_{1}},...,u_{w_{n}}}E\left(
L_{w_{1}}^{u_{w_{1}}}...L_{w_{n}}^{u_{w_{n}}}\right) ,$

\strut

for all $n>1$ in $\mathbb{N},$ where $d_{k}=\underset{v^{(k)}\in V(G:d_{k})}{%
\sum }q_{v^{(k)}}L_{v^{(k)}}\in D_{G}$ are arbitrary for $k=1,...,n..$
\end{theorem}

\strut

\begin{proof}
Suppose we have the compressed random variable $PaP$ of the $D_{G}$-valued
random variable $a.$ Then

\strut

$\ \ \ \ \ \ \ \ \ \ \ \ \ \ \ PaP=P_{V}(a)+P_{V}^{c}(a),$

\strut

where $P_{V}$ is the diagonal compression and $P_{V}^{c}$ is the
off-diagonal compression by $V\subset V(G).$ Clearly, we have the above
first cumulant of $PaP.$ Also, we have that

\strut

$\ k_{n}\left( d_{1}PaP,...,d_{n}PaP\right) $

\strut

$\ \ \ =k_{n}\left(
d_{1}(P_{V}(a)+P_{V}^{c}(a))~+...+~d_{n}(P_{V}(a)+P_{V}^{c}(a))\right) $

\strut

$\ \ \ =k_{n}\left( d_{1}P_{V}(a)~,...,~d_{n}P_{V}(a)\right) $

$\ \ \ \ \ \ \ \ \ \ \ \ \ \ \ \ \ \ \ \ \ \ \ +k_{n}\left(
d_{1}P_{V}^{c}(a)~,...,~d_{n}P_{V}^{c}(a)\right) $

\strut

by the $D_{G}$-freeness of $P_{V}(a)$ and $P_{V}^{c}(a)$

\strut

$\ \ \ =k_{n}\left( d_{1}P_{V}(a)~,...,~d_{n}P_{V}(a)\right) $

\strut

$\ \ \ \ \ \ \ \ \ \ \ +\underset{(i,j)\in \{1,...,N\}^{2},~i\neq j}{\sum }%
~k_{n}\left(
d_{1}(L_{v_{i}}aL_{v_{j}})~,...,~d_{n}(L_{v_{i}}aL_{v_{j}})\right) $

\strut

$\ \ \ =k_{n}\left( d_{1}P_{V}(a)~,...,~d_{n}P_{V}(a)\right) +0_{D_{G}}$

\strut

by Section 3.1

\strut

$\ \ \ =k_{n}\left( d_{1}P_{V}(a)~,...,~d_{n}P_{V}(a)\right) ,$

\strut

for all $n\in 2\mathbb{N}.$ Therefore, by [18], we can get the above result.
\end{proof}

\strut

\begin{remark}
The above theorem simply shows that

\strut

$\ \ \ \ \ \ k_{n}\left( d_{1}(PaP),...,d_{n}(PaP)\right) =k_{n}\left(
d_{1}P_{V}(a),...,d_{n}P_{V}(a)\right) ,$

\strut

for all $n\in \mathbb{N}$ and for any arbitrary $d_{1},...,d_{n}\in D_{G}.$
\end{remark}

\strut \strut

This says that the off-diagonal compressed part $P_{V}^{c}(a)$ does not
affect to compute the $D_{G}$-valued cumulants of the compressed random
variable $PaP.$ We can conclude that the compressed random variable $PaP$ of
the $D_{G}$-valued random variable by $P=\sum_{j=1}^{N}L_{v_{j}}$ and the
diagonal compressed random variable $P_{V}(a)$ of the random variable $a$ by 
$V=\{v_{1},...,v_{N}\}$ have the same distributions and hence they have the
same $D_{G}$-valued R-transforms.

\strut \strut \strut

\strut

\strut

\subsection{$D_{G}$-Freeness of Compressed Random Variables}

\strut

\strut \strut \strut

In this section, we will consider the $D_{G}$-freeness of compressed random
variables. In this section, we will consider the various conditions for the $%
D_{G}$-freeness of compressed random variables.

\strut

\begin{theorem}
Let $a,b\in \left( W^{*}(G),E\right) $ be a $D_{G}$-valued random variable
and let $V=\{v_{1},...,v_{N}\}$ be a finite subset of the vertex set $V(G)$
and $P=\sum_{j=1}^{N}L_{v_{j}},$ the corresponding projection in $W^{*}(G).$
Let $PaP$ and $PbP$ be the compressed random variable of $a$ and $b$ by $P$
in $\left( W^{*}(G),E\right) ,$ respectively. If $a$ and $b$ satisfy that

\strut

$\ \ \ \ \ \ \ \ \left( V\cap V(G:a)\right) \cap \left( V\cap V(G:b)\right)
=\emptyset ,$

\strut

$\ \ \ \ \ \ \ \ \ \ loop^{c}(G:a)\cap loop^{c}(G:b)=\emptyset ,\ \ \ $

and

$\ \ \ \ \ \ \ \ \ \ loop_{v_{i}}(G:a)\cap loop_{v_{j}}(G:b)=\emptyset ,$

\strut

for all choices $(i,j)\in \{1,...,N\}^{2},$ then $PaP$ and $PbP$ are free
over $D_{G}$ in $\left( W^{*}(G),E\right) .$
\end{theorem}

\strut

\begin{proof}
By the previous section, we have that

\strut

$\ \ \ \ \ \ PaP=P_{V}(a)+P_{V}^{c}(a)$ \ \ and \ \ $%
PbP=P_{V}(b)+P_{V}^{c}(b),$

\strut

where $P_{V},\,P_{V}^{c}:W^{*}(G)\rightarrow W^{*}(G)$ are the diagonal
compression and off-diagonal compression by $V,$ respectively. Moreover, $%
P_{V}(a)$ (resp. $P_{V}(b)$) and $P_{V}^{c}(a)$ (resp. $P_{V}^{c}(b)$) are
free over $D_{G}$ in $\left( W^{*}(G),E\right) .$ By the assumption and by
[18], $P_{V}(a)$ and $P_{V}(b)$ are free over $D_{G}$ in $\left(
W^{*}(G),E\right) .$ By the second condition $P_{V}^{c}(a)$ and $%
P_{V}^{c}(b) $ are free over $D_{G}$ in $\left( W^{*}(G),E\right) .$ This
shows that $\{P_{V}(a),P_{V}^{c}(a)\}$ and $\{P_{V}(b),P_{V}^{c}(b)\}$ are
free over $D_{G}$ in $\left( W^{*}(G),E\right) .$ So, $PaP$ and $PbP$ are
free over $D_{G}.$
\end{proof}

\strut \strut

\begin{theorem}
Let $a\in \left( W^{*}(G),E\right) $ be a $D_{G}$-valued random variable and
let $V_{1}=\{v_{1}^{(1)},...,v_{N_{1}}^{(1)}\}$ and $V_{2}=%
\{v_{1}^{(2)},...,v_{N_{2}}^{(2)}\}$ be finite subsets of the vertex set $%
V(G).$ Let $P$ and $Q$ be the corresponding projections of $V_{1}$ and $%
V_{2} $ in $W^{*}(G),$ respectively. Suppose that

\strut

$\ \ \ \ \ \ \ \left( V_{1}\cap V(G:a)\right) \cap \left( V_{2}\cap
V(G:a)\right) =\emptyset ,$

and

\ \ \ \ \ \ \ \ $loop_{v_{i}^{(1)}}(G:a)\cap
loop_{v_{j}^{(2)}}(G:a)=\emptyset ,$

\strut \strut

for all choices $(i,j)\in \{1,...,N_{1}\}\times \{1,...,N_{2}\}.$ Then the
corresponding diagonal compressed random variables $PaP$ and $QaQ$ of $a$ by 
$P$ and $Q$ are free over $D_{G}$ in $\left( W^{*}(G),E\right) .$
\end{theorem}

\strut

\begin{proof}
By hypothesis and by [18], $P_{V_{1}}(a)$ and $Q_{V_{2}}(a)$ are free over $%
D_{G}$ in $\left( W^{*}(G),E\right) ,$ where $P_{V_{1}},\,%
\,Q_{V_{2}}:W^{*}(G)\rightarrow W^{*}(G)$ are diagonal compressions by $%
V_{1} $ and $V_{2},$ respectively.
\end{proof}

\strut \strut

\strut

\strut

\begin{quote}
\textbf{Reference}

\strut

\strut

{\small [1] \ \ A. Nica, R-transform in Free Probability, IHP course note,
available at www.math.uwaterloo.ca/\symbol{126}anica.}

{\small [2]\strut \ \ \ A. Nica and R. Speicher, R-diagonal Pair-A Common
Approach to Haar Unitaries and Circular Elements, (1995), www
.mast.queensu.ca/\symbol{126}speicher.\strut }

{\small [3] \ }$\ ${\small B. Solel, You can see the arrows in a Quiver
Operator Algebras, (2000), preprint}

{\small \strut [4] \ \ A. Nica, D. Shlyakhtenko and R. Speicher, R-cyclic
Families of Matrices in Free Probability, J. of Funct Anal, 188 (2002),
227-271.}

{\small [5] \ \ D. Shlyakhtenko, Some Applications of Freeness with
Amalgamation, J. Reine Angew. Math, 500 (1998), 191-212.\strut }

{\small [6] \ \ D.Voiculescu, K. Dykemma and A. Nica, Free Random Variables,
CRM Monograph Series Vol 1 (1992).\strut }

{\small [7] \ \ D. Voiculescu, Operations on Certain Non-commuting
Operator-Valued Random Variables, Ast\'{e}risque, 232 (1995), 243-275.\strut 
}

{\small [8] \ \ D. Shlyakhtenko, A-Valued Semicircular Systems, J. of Funct
Anal, 166 (1999), 1-47.\strut }

{\small [9] \ \ D.W. Kribs and M.T. Jury, Ideal Structure in Free
Semigroupoid Algebras from Directed Graphs, preprint}

{\small [10]\ D.W. Kribs and S.C. Power, Free Semigroupoid Algebras, preprint%
}

{\small [11]\ I. Cho, Amalgamated Boxed Convolution and Amalgamated
R-transform Theory, (2002), preprint.}

{\small [12] I. Cho, The Tower of Amalgamated Noncommutative Probability
Spaces, (2002), Preprint.}

{\small [13]\ I. Cho, Compatibility of a Noncommutative Probability Space
and a Noncommutative Probability Space with Amalgamation, (2003), Preprint}

{\small [14] I. Cho, An Example of Scalar-Valued Moments, Under
Compatibility, (2003), Preprint.}

{\small [15] I. Cho, Free Semigroupoid Probability Theory, (2004), Preprint.}

{\small [16] I. Cho, Graph }$W^{*}$-{\small Probability Theory, (2004),
Preprint.}

{\small [17] I. Cho, Random Variables in Graph }$W^{*}$-{\small Probability
Spaces, (2004), Preprint.}

{\small [18] I. Cho, Diagonal Compressed Random Variables in the Graph }$%
W^{*}$-{\small Probability Space, (2004), Preprint.}

{\small [19] I. Cho, Free Product of Two Graph }$W^{*}${\small -Probability
Spaces, (2004), Preprint.}

{\small [20] P.\'{S}niady and R.Speicher, Continous Family of Invariant
Subspaces for R-diagonal Operators, Invent Math, 146, (2001) 329-363.}

{\small [21] R. Speicher, Combinatorial Theory of the Free Product with
Amalgamation and Operator-Valued Free Probability Theory, AMS Mem, Vol 132 ,
Num 627 , (1998).}

{\small [22] R. Speicher, Combinatorics of Free Probability Theory IHP
course note, available at www.mast.queensu.ca/\symbol{126}speicher.\strut }
\end{quote}

\end{document}